\documentclass[a4paper]{article}
\usepackage{amsfonts}
\usepackage{amssymb}
\usepackage{amsmath}
\usepackage{latexsym}
\usepackage{graphicx}
\usepackage{amsthm}
\usepackage{calrsfs}
\newtheorem{theorem}{Theorem}[section]

\newtheorem{problem}[theorem]{Problem}
\newtheorem{proposition}[theorem]{Proposition}
\newtheorem{example}[theorem]{Example}
\newtheorem{definition}[theorem]{Definition}
\newtheorem{remark}[theorem]{Remark}

\begin{document}

\parbox{1mm}

\begin{center}
{\bf {\sc \Large Probabilistic-valued decomposable set functions
\\ with respect to triangle functions}}
\end{center}

\vskip 12pt

\begin{center}
{\bf Lenka HAL\v CINOV\'A, Ondrej HUTN\'IK and Jana
MOLN\'AROV\'A}
\end{center}


\hspace{5mm}\parbox[t]{12cm}{\fontsize{9pt}{0.1in}\selectfont\noindent{\bf
Abstract.} In the framework of generalized measure theory the
probabilistic-valued decomposable set functions are introduced
with triangle functions $\tau$ in an appropriate probabilistic
metric space as natural candidates for the "addition", leading to
the concept of $\tau$-decomposable measures. Several set
functions, among them the classical (sub)measures, previously
defined $\tau_T$-submeasures~\cite{HutMes},
$\tau_{L,A}$-submeasures~\cite{HalHutMes2} as well as recently
introduced Shen's (sub)measures~\cite{Shen} are described and
investigated in a unified way. Basic properties and
characterizations of $\tau$-decomposable (sub)measures are also
studied and numerous extensions of results from the above
mentioned papers are provided. } \vskip 24pt

\section{Introduction}


Real world applications often require dealing with such situations
when the exact numerical values of (sub)measure of a set may not
be provided, but at least some probabilistic assignment still
could be done. A similar situation is discussed in the framework
of information measures in~\cite{KdF74}. For instance, consider a
grant agency providing a financial support for research in some
area. From the set of all grant applications only "successful"
(depending on some internal rules of agency) will receive certain
amount of money. So, we have only a probabilistic information
about measure of the set of "successful" grant applications. Of
course, the knowledge of this information depends on many
different aspects: total budget of money to be divided, internal
rules of agency, quality of reviewers (if any), etc. Further
examples are provided by lotteries, or guessing results when we
have a probabilistic information about (counting) measure of
possibilities to win the prize. A closely related concept can be
found in Moore's interval mathematics~\cite{MB}, where the use of
intervals in data processing is due to measurement inaccuracy and
due to rounding. Here intervals can be considered in distribution
function form linked to random variables uniformly distributed
over the relevant intervals. These model examples resemble the
original idea of Menger of PM-spaces, see~\cite{Menger}, where the
replacement of a positive number by a distance distribution
function was motivated by thinking of situations where the exact
distance between two objects may not be provided, but some
probability assignment is still possible. Thus, the
importance/diameter/measure of a set might be represented by a
distance distribution function. Recently, probabilistic approaches
were successfully applied to modelling uncertain logical
arguments~\cite{Hunter}, to approximations of incomplete
data~\cite{GCK}, to inference rules playing an important role in
non-monotonic reasoning~\cite{GS}, or to cluster structure
ensemble~\cite{cinania}.

In the paper by Hutn\'ik and Mesiar~\cite{HutMes} the notion of
\textit{$\tau_T$-submeasure} was defined intended to be a certain
(non-additive) set function $\gamma$ on a ring $\Sigma$ of subsets
of a non-empty set $\Omega$ taking values in the set $\Delta^+$ of
distribution functions of non-negative random variables satisfying
$\gamma_\emptyset=\varepsilon_0$ , "antimonotonicity" property
$\gamma_{E}\geq \gamma_F$ whenever $E,F\in\Sigma$ with $E\subseteq
F$, and "subadditivity" property of the form
\begin{equation}\label{subadditivityT}
\gamma_{E\cup F}(x+y) \geq T(\gamma_E(x), \gamma_F(y)), \quad
E,F\in\Sigma, x,y>0,
\end{equation} with $T$ being a
left-continuous t-norm. Here, $\varepsilon_0$ is the distribution
function of Dirac random variable concentrated at point 0. As it
is shown in~\cite{HutMes}, such $\tau_T$-submeasures can be seen
as fuzzy number-valued submeasures. In this case the value
$\gamma_E$ may be seen as a non-negative LT-fuzzy number,
see~\cite{DKMP}, where $\tau_T(\gamma_E,\gamma_F)$ corresponds to
the $T$-sum of fuzzy numbers $\gamma_E$ and $\gamma_F$. Also, each
$\tau_M$-submeasure $\gamma$ with the minimum t-norm $M(x,y) =
\min\{x,y\}$ can be represented by means of a non-decreasing
system $(\eta_\alpha)_{\alpha \in [0,1]}$ of numerical submeasures
as follows
$$\gamma_E(x) = \sup\bigl\{\alpha \in [0,1];\, \eta_\alpha(E)\leq x\bigr\}, \quad E\in\Sigma.$$
This representation resembles the horizontal representation
$(S_\alpha)_{\alpha \in [0,1]}$ of a fuzzy subset $S$.

The study of probabilistic-valued set functions continued in
papers~\cite{HalHutMes2} and~\cite{HalHutMes}, where a more
general concept has been used. In fact, in~\cite{HalHutMes} a
generalization of $\tau_{T}$-sub\-measu\-res was suggested which
involves suitable operations $L$ replacing the standard addition
$+$ on $\overline{\mathbb{R}}_+$. On the other hand, since t-norms
are rather special operations on the unit interval $[0, 1]$, the
paper~\cite{HalHutMes2} deals with a number of possible
generalizations based on aggregation functions in general studying
certain properties of the corresponding probabilistic
(sub)measures and their (sub)measure spaces.

Recently, in~\cite{Shen} Shen defined and studied a class of
probabilistic (sub)measures which seem they do not fit to the
concept of previously mentioned results of Hutn\'ik and Mesiar. In
fact, Shen's definition of a probabilistic-valued
$\top$-decomposable supmeasure\footnote{\textit{supmeasure} in the
terminology of Shen corresponds to \textit{submeasure} in our
terminology, see~\cite[Definition 4.1(v)]{Shen}} $\mathfrak{M}:
\Sigma\to\Delta^+$ with the "subadditivity" property
\begin{equation}\label{subadditivityShen}
\mathfrak{M}_{E\cup
F}(t)\geq\top(\mathfrak{M}_{E}(t),\mathfrak{M}_{F}(t)), \quad
E,F\in\Sigma, t>0,
\end{equation} with $\top$ being a t-norm,
corresponds to the notion of $\tau_{\max,T}$-submeasure defined
in~\cite{HalHutMes}, however a deeper contextual understanding was
still unclear from that paper. Thus, in this paper we provide a
deep insight into all the mentioned notions being, in fact,
special cases of a \textit{probabilistic-valued set function with
respect to a triangle function}. Recall that a triangle function
$\tau$ is a binary operation on $\Delta^+$ such that the triple
$(\Delta^+, \tau, \leq)$ forms a commutative, partially ordered
semigroup with neutral element $\varepsilon_0$.

More precisely, the notion of $\tau_T$-submeasure is related to
the triangle function $\tau=\tau_T$ given by
\begin{equation}\label{tauT}
\tau_T(G,H)(x)=\sup_{u+v=x} T(G(u),H(v)), \quad G,H\in\Delta^+
\end{equation} with $T$ being a left-continuous
t-norm. Thus, the "subadditivity" property~(\ref{subadditivityT})
resembles the "probabilistic analogue" of the triangle inequality
in the Menger probabilistic metric space (under $T$),
see~\cite{SS}. Moreover, $\tau_{L,T}$-submeasures defined
in~\cite[Definition 2.3]{HalHutMes} are related to the (triangle)
function
\begin{equation}\label{tauLT}
\tau_{L,T}(G,H)(x) = \sup_{L(u,v)=x} T(G(u),H(v)), \quad
G,H\in\Delta^+,
\end{equation} with a suitable operation $L$ on
$\overline{\mathbb{R}}_+$. Even more, Shen's considerations are
related to the pointwisely defined (triangle) function
$$\tau_\top(G,H)(t)=\top(G(t), H(t)), \quad
G,H\in\Delta^+,$$ and the "subadditivity"
property~(\ref{subadditivityShen}) is related to the triangle
inequality of the corresponding probabilistic metric space. So, we
can see that triangle functions are the main ingredient which
connects all the mentioned notions of (sub)measure. Thus,
considering a general triangle function $\tau$ on $\Delta^+$ we
define and study certain properties of $\tau$-decomposable
(sub)measures on a ring $\Sigma$ of subsets of
$\Omega\neq\emptyset$ in this general setting.

In the next section the short overview of basic notions and
definitions is given. In Section~\ref{secsetfunctions} we
introduce the basic object of our study: a $\tau$-decomposable set
function with values in distance distribution functions and
provide a number of concrete examples. Several properties of such
set functions are then investigated in Section~\ref{secproperties}
and results related to probabilistic Hausdorff distance are
provided in Section~\ref{sec_hausdorff} generalizing the recent
results of Shen~\cite{Shen}.

\section{Basic notions and definitions}\label{sectionprem}

In order to make the exposition self-contained, here we remind the
reader the basic notions and constructions used in this paper.

\paragraph{Distribution functions} Let $\Delta$ be the family of all distribution functions on the
extended real line $\overline{\mathbb{R}}:=[-\infty,+\infty]$,
i.e., $F: \overline{\mathbb{R}} \to [0,1]$ is non-decreasing, left
continuous on the real line $\mathbb{R}$ with $F(-\infty)=0$ and
$F(+\infty)=1$. A \textit{distance distribution function} is a
distribution function whose support is a subset of
$\overline{\mathbb{R}}_+:=[0,+\infty]$, i.e., a distribution
function $F: \overline{\mathbb{R}}\to [0,1]$ with $F(0) = 0$. The
class of all distance distribution functions will be denoted by
$\Delta^+$.

For a distance distribution function $G$ and a non-negative
constant $c\in \overline{\mathbb{R}}_{+}$ we define the
multiplication of $G$ by a constant $c$ as follows
\begin{equation}\label{operaciestf1}
\left(c\odot G\right)(x):=
\begin{cases} \varepsilon_0(x), & c\in \{0,+\infty\}, \\
G\left(\frac{x}{c}\right), & \textrm{otherwise}.
\end{cases}\end{equation} \noindent Clearly, $c\odot G\in\Delta^+$.

A \textit{triangle function} is a mapping $\tau: \Delta^+ \times
\Delta^+ \to \Delta^+$ which is symmetric, associative,
non-decreasing in each variable and has $\varepsilon_0$ as the
identity, where $\varepsilon_0$ is the distribution function of
Dirac random variable concentrated at point 0. More precisely, for
$a\in[-\infty,+\infty[$ we put
$$\varepsilon_a(x) :=
\begin{cases}
1& \textrm{for}\,\, x>a, \\
0 & \textrm{otherwise}.
\end{cases}$$ For more
details on triangle functions we recommend the overview
paper~\cite{SamSem}. The most important triangle functions are
those obtained from certain aggregation functions, especially
t-norms.

\paragraph{Aggregation functions}
A binary \textit{aggregation function} $A: [0,1]^2\to [0,1]$ is a
non-decreasing function in both components with the boundary
conditions $A(0,0)=0$ and $A(1,1)=1$. The class of all binary
aggregation functions will be denoted by $\mathcal{A}$. For more
details on aggregation functions we recommend the
monograph~\cite{GPMM}.

A \textit{triangular norm}, shortly a t-norm, is a commutative
lattice ordered semi-group on $[0,1]$ with identity 1. The most
important are the minimum t-norm $M(x,y) := \min\{x,y\}$, the
product t-norm $\Pi(x,y) := xy$, the {\L}ukasiewicz t-norm $W(x,y)
:= \max\{x+y-1,0\}$, and the drastic product t-norm
$$D(x,y) :=
\begin{cases}
\min\{x,y\} & \textrm{for}\,\, \max\{x,y\}=1 \\
0 & \textrm{otherwise}.
\end{cases}$$ For more information about t-norms and their properties
we refer to books~\cite{KMP,SS}. Throughout this paper
$\mathcal{T}$ denotes the class of all t-norms.

\paragraph{Binary operations on non-negative reals} Let us denote by $\mathcal{L}$ the set of all
binary operations $L$ on $\overline{\mathbb{R}}_+$ such that

\begin{itemize} \item[(i)] $L$ is commutative and associative; \item[(ii)] $L$ is jointly strictly
increasing, i.e., for all $u_1, u_2, v_1,
v_2\in\overline{\mathbb{R}}_+$ with $u_1<u_2$, $v_1<v_2$ holds
$L(u_1,v_1)<L(u_2,v_2)$; \item[(iii)] $L$ is continuous on
$\overline{\mathbb{R}}_+\times\overline{\mathbb{R}}_+$;
\item[(iv)] $L$ has $0$ as its neutral element.
\end{itemize}\noindent Observe that $L\in\mathcal{L}$
is a jointly increasing pseudo-addition on
$\overline{\mathbb{R}}_+$ in the sense of~\cite{SM}. The usual
(class of) examples of operations in $\mathcal{L}$ are
\begin{align*}
K_\alpha(x,y) & := (x^\alpha+y^\alpha)^{\frac{1}{\alpha}}, \quad
\alpha>0, \\ K_\infty(x,y) & := \max\{x,y\}.
\end{align*}In general, $L\in\mathcal{L}$ if and only if there is a (possibly
empty) system $(]a_k,b_k[)_{k \in K}$ of pairwise disjoint open
subintervals of $]0,+\infty[$, and a system $(\ell_k)_{k \in K}$
of increasing bijections $\ell_k: [a_k,b_k]\to
\overline{\mathbb{R}}_+$ so that
$$L(x,y) =
\begin{cases}
\ell_{k}^{-1}(\ell_k(x)+\ell_k(y)) & \mbox{if } (x,y) \in ]a_k,b_k[^2,\\
\max\{x,y\} & \mbox{otherwise}.
\end{cases}$$ For more details see~\cite{KMP}.
Further, for $(L,A)\in\mathcal{L}\times\mathcal{A}$ define the
function
\begin{equation}\label{tauLA}
\tau_{L,A}(G,H)(x) := \sup_{L(u,v)=x} A(G(u), H(v)), \quad G,H\in
\Delta^+.\end{equation} However, $\tau_{L,A}$ need not be
associative in general, thus need not be a triangle function, but
it has good properties on $\Delta^+$. As it is shown
in~\cite{SamSem}, the appropriate choice for $A$ is a semi-copula
$S$ -- a binary aggregation function on $[0,1]$ with 1 as its
neutral element. Indeed, the left-continuity of $S$ guarantees
that $\tau_{L,S}$ is a binary operation on $\Delta^+$,
cf.~\cite[Lemma 7.1]{SamSem}. Also, for any semi-copula $S$,
$$\tau_{K_\infty,S}(G,H)(x) = S(G(x),H(x))$$ is the operation pointwisely induced by $S$ on $\Delta^+$.
For more information about (triangular) functions in connection
with various aggregation functions and their properties we refer
to~\cite{SamSem}.

\paragraph{Probabilistic metric spaces} The
function~(\ref{tauLA}), at least for a special choice of $L$ and
$A$, naturally arises in the context of probabilistic metric
spaces. Recall that a \textit{probabilistic metric space}
(PM-space, for short) is a non-empty set $\Omega$ together with a
family $\mathcal{F}$ of probability functions $F_{p,q}(x)$
(interpreted as the probability that distance between elements
$p,q$ of $\Omega$ is less than $x$) with $F_{p,q}(0)=0$ satisfying
\begin{itemize}
\item[(i)] $F_{p,q} = \varepsilon_0$ if and only if $p = q$;
\item[(ii)] $F_{p,q} = F_{q,p}$; \item[(iii)] $F_{p,r} \geq
\tau(F_{p,q}, F_{q,r})$,
\end{itemize}\noindent where $\tau$ is a triangle
function on $\Delta^+$. A particular case includes the triangle
inequality
\begin{equation}\label{triangleineq}
F_{p,r}(L(x,y)) \geq A(F_{p,q}(x), F_{q,r}(y))
\end{equation} which
holds for all $p,q,r \in \Omega$ and all real $x,y$ with a
suitable binary aggregation function $A$ and operation $L$. In
particular, the triple $(\Omega, \mathcal{F}, \tau_{L,T})$ with
$\tau_{L,T}$ given by~(\ref{tauLA}) for $A=T$ (a left-continuous
t-norm) is called an \textit{$L$-Menger PM-space} (under $T$). For
$L=K_\infty$ we get a non-Archimedean Menger PM-space $(\Omega,
\mathcal{F}, \tau_{K_\infty, T})$ with $T\in\mathcal{T}$. Note
that if $\mathcal{F}$ satisfies $F_{p,p}=\varepsilon_0$ for each
$p\in\Omega$ and properties~(ii) and~(iii), then the triple
$(\Omega,\mathcal{F}, \tau)$ will be called a
\textit{probabilistic
pseudo-metric space} (PpM-space, for short). 
In general, different triangle functions lead to PM-spaces with
different geometric and topological properties.

\section{Probabilistic-valued decomposable set functions w.r.t. a triangle
function}\label{secsetfunctions}

As we have already mentioned in the introduction, a natural origin
of probabilistic-valued set functions comes from the fact that
they work in such situations in which we have only a
\textit{probabilistic information} about measure of a set. For
example, if rounding of reals is considered, then the uniform
distributions over intervals describe our information about the
measure of a set. Thus, the probabilistic-valued (sub)measures
represent the concept of (sub)measures probabilistically rather
than deterministically.

Here we introduce the basic notion of probabilistic decomposable
(sub)measure in its general form. For better readability we also
use the following conventions:
\begin{itemize}
\item[(i)] for a probabilistic-valued set function $\gamma: \Sigma
\to \Delta^+$ we write $\gamma_E(x)$ instead of $\gamma(E)(x)$;
\item[(ii)] since $\Delta^+$ is the set of all distribution
functions with support $\overline{\mathbb{R}}_+$, we state the
expression for a mapping $\gamma: \Sigma \to \Delta^+$ just for
positive values of $x$. In case $x\leq 0$ we always suppose
$\gamma_{\cdot}(x)=0$.
\end{itemize}

\begin{definition}\rm\label{defLA-sub}
Let $\tau$ be a triangle function on $\Delta^+$ and $\Sigma$ be a
ring of subsets of $\Omega\neq \emptyset$. A mapping $\gamma:
\Sigma \to \Delta^+$ with $\gamma_\emptyset = \varepsilon_0$ is
said to be a $\tau$-\textit{decomposable submeasure}, if
$\gamma_{E\cup F} \geq \tau(\gamma_E, \gamma_F)$ for each disjoint
sets $E,F \in \Sigma$. If in the preceding inequality equality
holds, then $\gamma$ is said to be a $\tau$-\textit{decomposable
measure} on $\Sigma$.
\end{definition}

\begin{remark}\rm
In fact, a mapping $\gamma: \Sigma\to \Delta^+$ is a
probabilistic-valued set function, where the value of $\gamma_E$
at $x$ may be interpreted as the probability that a numerical
(sub)measure of the set $E$ is less than $x$. A triangle function
$\tau$ is a natural choice for "aggregation" of $\gamma_E$ and
$\gamma_F$ in order to compare them with $\gamma_{E\cup F}$. For a
$\tau$-decomposable measure $\gamma$ we naturally expect that
$\gamma_{E\cup F}$ is the same distance distribution function as
$\gamma_{F\cup E}$ for disjoint sets $E,F\in\Sigma$, from which
follows that $\tau(\gamma_E,\gamma_F)=\tau(\gamma_F,\gamma_E)$,
i.e., $\tau$ has to be commutative. Moreover, from the natural
equality $\gamma_{(E\cup F)\cup G} = \gamma_{E\cup(F\cup G)}$ we
obtain
$\tau(\tau(\gamma_E,\gamma_F),\gamma_G)=\tau(\gamma_E,\tau(\gamma_F,\gamma_G))$,
i.e., $\tau$ has to be associative. Since
$\gamma_E=\gamma_{E\cup\emptyset}=\tau(\gamma_E,\gamma_\emptyset)=\tau(\gamma_E,\varepsilon_0)$,
then $\varepsilon_0$ has to be neutral element of $\tau$. The role
of monotonicity of $\tau$ (a non-decreasing function in each
place) will be examined in what follows.
\end{remark}

\begin{theorem}
Let $\tau$ be a triangle function on $\Delta^+$ and $\Sigma$ be a
ring of subsets of $\Omega\neq \emptyset$. Then each
$\tau$-decomposable measure $\gamma$ is "antimonotone" on
$\Sigma$, i.e., $\gamma_E\geq\gamma_F$ whenever $E,F\in\Sigma$
such that $E\subseteq F$.
\end{theorem}

\proof From monotonicity of triangular function $\tau$ with
$\varepsilon_0$ as identity we easily have
$$\gamma_F = \gamma_{E\cup (F\setminus E)} = \tau(\gamma_{E},
\gamma_{F\setminus E}) \leq \tau(\gamma_{E},\varepsilon_0) =
\gamma_{E}$$ whenever $E, F\in\Sigma, E\subseteq F$.\qed \vskip
5pt


\begin{remark}\rm
The probabilistic interpretation of this property is as follows:
the probability that a numerical measure of the set $E$ is less
than $x$ is greater than or equal to the probability that a
numerical measure of the set $F$ is less than $x$.
\end{remark}

"Antimonotonicity" property does not hold for arbitrary
$\tau$-decomposable submeasures $\gamma: \Sigma \to \Delta^+$.
Therefore, in~\cite[Definition 2.1]{HalHutMes2} we have considered
the notion of probabilistic submeasure w.r.t. a function
$\tau=\tau_{L,A}$ being a $\tau_{L,A}$-decomposable antimonotone
submeasure on $\Sigma$ where the "subadditivity" property
$\gamma_{E\cup F} \geq \tau_{L,A}(\gamma_E, \gamma_F)$ holds for
arbitrary sets $E,F\in\Sigma$. Therefore we state the following
easy observation.

\begin{theorem}\label{thmantimonotone}
Let $\tau$ be a triangle function on $\Delta^+$ and $\Sigma$ be a
ring of subsets of $\Omega\neq \emptyset$. If $\gamma: \Sigma \to
\Delta^+$ is a $\tau$-decomposable antimonotone submeasure, then
the inequality $\gamma_{E\cup F} \geq \tau(\gamma_E, \gamma_F)$
holds for arbitrary sets $E,F \in \Sigma$.
\end{theorem}

\proof For $E\cap F=\emptyset$ the inequality holds by
Definition~\ref{defLA-sub}. Let $E,F \in \Sigma$ and $E\cap
F\neq\emptyset$. Then the inequality $$\gamma_{E\cup F} =
\gamma_{E\cup ((E\cup F)\setminus E)}\geq \tau(\gamma_E,
\gamma_{(E\cup F)\setminus E})\geq \tau(\gamma_E, \gamma_F)$$
follows from antimonotonicity of $\gamma$ and monotonicity of
triangle functions. \qed \vskip 5pt

Now we are in position to provide some examples of
$\tau$-decomposable measures based on various constructions and
related to different classes of aggregation functions appearing in
the definition of the underlying triangle function. Note that
under the "numerical (sub)measure" we mean a real-valued set
function with the (usual) (sub)additivity property.

\begin{example}\rm
Let $m$ be a numerical (additive) measure on $\Sigma$. If $m$ is
$L$-decomposable, i.e., $m(E\cup F)=L(m(E),m(F))$ with
$L\in\mathcal{L}$ and disjoint sets $E,F\in\Sigma$, then for any
$\Phi\in\Delta^+$ the set function $\gamma^\Phi:\Sigma\to\Delta^+$
defined by $\gamma^\Phi_E:=m(E)\odot \Phi$ is a
$\tau_{L,M}$-decomposable measure on $\Sigma$, where $\tau_{L,M}$
is given by~(\ref{tauLT}) with the minimum t-norm $T=M$. Indeed,
"additivity" property follows from~\cite[Section 7.7]{SS}: the
function $\tau_{L,M}$ is the only triangular function with the
property
$$\tau_{L,M}\left(c_1\odot H, c_2\odot H\right)= L(c_1,c_2)\odot
H,\qquad c_{1}, c_{2}\in\overline{\mathbb{R}}_{+}, H\in
\Delta^+.$$ Notice that $\gamma^\Phi$ need not be a
$\tau_{L,T}$-decomposable measure for $T\neq M$.
\end{example}


\begin{example}\rm
\textit{Shen's $\top$-probabilistic decomposable measures},
cf.~\cite{Shen}: this class of measures
$\mathfrak{M}:\Sigma\to\Delta^+$ of the form $\mathfrak{M}_{E\cup
F}(t)=\top(\mathfrak{M}_{E}(t),\mathfrak{M}_{F}(t))$ for disjoint
$E,F\in\Sigma$ corresponds to the class of $\tau$-decomposable
measures w.r.t. the triangle function $\Pi_\top: \Delta^{+}\times
\Delta^{+}\rightarrow \Delta^{+}$ of the form
\begin{equation}\label{tautop}
\Pi_\top(G,H)(t)=\top(G(t), H(t)), \quad
G,H\in\Delta^+,\end{equation} with $\top$ being a t-norm, or,
equivalently, to a $\tau_{\max,T}$-decomposable measure. According
to~\cite[Theorem 5.2]{SamSem} left-continuity of t-norm $\top$ is
a necessary and sufficient condition for $\Pi_\top$ being a
triangle function. Note that Shen does not consider the
left-continuity of $\top$ in his definition of a
$\top$-probabilistic decomposable measure, cf.~\cite[Definition
4.1]{Shen}. Moreover, Shen's considerations are made on a
$\sigma$-algebra instead of a ring of subsets of
$\Omega\neq\emptyset$ to provide a countable extension of
$\top$-decomposable measures.
\end{example}


\begin{example}\rm
Convolution of distance distribution functions is a binary
operation $*$ on $\Delta^+$ given by
$$\tau_*(G,H)(x):= (G*H)(x)=
\begin{cases}
0, & x=0, \\
\int_{0}^x G(x-t)\,\textrm{d}H(t), & x\in]0,+\infty[, \\
1, & x=+\infty,
\end{cases}$$ for each $G,H\in\Delta^+$, where the integral is
meant in the sense of Lebesgue-Stieltjes. According
to~\cite[Theorem 13.4 and Definition 13.1]{SamSem} $\tau_*$ is a
triangle function on $\Delta^+$. The corresponding
$\tau_*$-decomposable (sub)measure $\gamma$ provides an extension
of a notion of (sub)measure to (sub)measures which can be used in
the Wald spaces -- those involving the convolution $*$ of distance
distribution functions, but also in a wider class of PM-spaces.

If $m$ is a numerical measure on $\Sigma$, then the set function
$\gamma:\Sigma\to\Delta^+$ given by $\gamma_E=\varepsilon_{m(E)}$
for $E\in\Sigma$ is a $\tau_*$-decomposable measure, because
$$\gamma_{E\cup F} = \varepsilon_{m(E\cup F)} =
\varepsilon_{m(E)+m(F)} =
\tau_*(\varepsilon_{m(E)},\varepsilon_{m(F)}) =
\tau_*(\gamma_E,\gamma_F),$$ cf.~\cite[Theorem 14.1(f)]{SamSem}.
\end{example}


The latter example of mapping $\gamma_E=\varepsilon_{m(E)}$ is
also interesting from the viewpoint of other triangle functions.
It gives rise claiming that \textit{each numerical (additive)
measure can be regarded as a probabilistic-valued decomposable
measure}. Thus, $\tau$-decomposable measures are extensions of the
classical measure.

\begin{example}\rm
Let $\mu: \Sigma\to[0,+\infty)$ be a set function with
$\mu(\emptyset)=0$ and for $E\in\Sigma$ put
$\gamma_E=\varepsilon_{\mu(E)}$. Due to~\cite[Theorem
14.1]{SamSem} we have the following results:
\begin{itemize}
\item[(i)] if $\mu$ is additive, then $\gamma$ is a
$\tau_D$-decomposable measure; \item[(ii)] if $\mu$ is additive
and $T\in\mathcal{T}$ is continuous, then $\gamma$ is a
$\tau_T$-decomposable measure; \item[(iii)] if $\mu$ is
$L$-decomposable with $L\in\mathcal{L}$, and $T\in\mathcal{T}$ is
continuous, then $\gamma$ is a $\tau_{L,T}$-decomposable measure;
\item[(iv)] if $\mu$ is $L$-decomposable with $L\in\mathcal{L}$
and $Q$ is a symmetric quasi-copula with an associative dual
quasi-copula $\overline{Q}$, then $\gamma$ is a
$\rho_{L,Q}$-decomposable measure, where
$$\rho_{L,Q}(G,H)(x) = \inf_{L(u,v)=x} \overline{Q}(G(u),H(v)),
\quad G,H\in\Delta^+;$$ 
\item[(v)] if $\mu$ is $K_\infty$-decomposable, then $\gamma$ is
$\Pi_\top$-decomposable measure for each $\top\in\mathcal{T}$,
where $\Pi_\top$ is given by~(\ref{tautop}), i.e., $\gamma$ is a
Shen's $\top$-probabilistic decomposable measure. Note that if
$\mu$ is additive, then $\gamma$ need not be a
$\Pi_\top$-decomposable measure. For instance, considering
$\top=M$ (the minimum t-norm) yields
$$\varepsilon_{\mu(E)+\mu(F)}(x) \neq
M(\varepsilon_{\mu(E)}(x),\varepsilon_{\mu(F)}(x))$$ whenever
$\mu(E)\neq \mu(F)\neq\emptyset$.
\end{itemize}
\end{example}


\section{Properties and constructions of decomposable (sub)measures}\label{secproperties}

Clearly, each $\tau$-decomposable measure on $\Sigma$ is a
$\tau$-decomposable submeasure on $\Sigma$. The following result
provides a characterization of the class of $\tau$-decomposable
measures on $\Sigma$. In fact, it provides a generalization
of~\cite[Theorem 4.1]{Shen} to an arbitrary triangle function.

\begin{theorem}
Let $\tau$ be a triangle function on $\Delta^+$. Then $\gamma$ is
a $\tau$-decomposable measure on $\Sigma$ if and only if
$\tau(\gamma_{E\cup F}, \gamma_{E\cap F})=\tau(\gamma_E,\gamma_F)$
for each $E,F\in\Sigma$.
\end{theorem}

\proof "$\Leftarrow$" If $E\cap F=\emptyset$, then $\gamma_{E\cap
F}=\varepsilon_0$ and
$$\tau(\gamma_E, \gamma_F)=\tau(\gamma_{E\cup F}, \gamma_{E\cap F}) =
\tau(\gamma_{E\cup F}, \varepsilon_0) = \gamma_{E\cup F}.$$

\noindent "$\Rightarrow$" Since $E\cup F = (E\cap
F)\cup(E\setminus F)\cup(F\setminus E)$, then we have
\begin{align*}
\tau(\gamma_{E\cap F},\gamma_{E\cup F}) &= \tau(\gamma_{E\cap
F},\tau(\gamma_{E\cap F},\tau(\gamma_{E\setminus
F},\gamma_{F\setminus E})) = \tau(\gamma_{E\cap
F},\tau(\tau(\gamma_{E\cap F},\gamma_{E\setminus
F}),\gamma_{F\setminus E})) \\ & = \tau(\gamma_{E\cap
F},\tau(\gamma_{E},\gamma_{F\setminus E})) = \tau(\gamma_{E\cap
F},\tau(\gamma_{E},\gamma_{F\setminus E})) =
\tau(\gamma_{E},\tau(\gamma_{E\cap F},\gamma_{F\setminus E})) \\
& = \tau(\gamma_E,\gamma_F)
\end{align*}which completes the proof.
\qed \vskip 5pt


Now we are interested in a question to generate new decomposable
(sub)measures from given ones. For that reason we will need the
following notion, cf.~\cite[Definition 6.1]{SamSemII}: Let
$(X,\leq)$ be a partially ordered set and let $f$ and $g$ be two
binary operations on $X$. Then \textit{$f$ dominates $g$}, written
$f \gg g$, if, for all $x, y, u, v \in X$,
$$f(g(x, y), g(u, v)) \geq g(f(x, u), f(y, v)).$$ It is well known that dominance is an antisymmetric and reflexive,
but not transitive relation on the set of triangle functions,
cf.~\cite[Corollary 6.1]{SamSemII}. Note that the idea of
dominance appears in much previous paper~\cite{SMB} in the context
of aggregation functions, or even in~\cite{dBM} in the context of
triangular norms.

A triangle function $\tau$, such that for each
$c\in\overline{\mathbb{R}}_{+}$ and each $G,H\in\Delta^+$ it holds
$$c\odot\tau(G,H) = \tau(c\odot G, c\odot H)$$ will
be called a \textit{distributive triangle function}. Recall that
the operation $\odot$ is defined by~(\ref{operaciestf1}). The
standard examples of distributive triangle functions are $\tau_T$
and $\Pi_\top$.

\begin{theorem}
Let $\tau, \vartheta$ be two triangle functions on $\Delta^+$ and
$\gamma^1, \gamma^2:\Sigma\to\Delta^+$ be $\tau$-decomposable
measures. Then
\begin{itemize}
\item[(i)] if $\tau$ is distributive, the set function
$\gamma:=c\odot\gamma^1$ is a $\tau$-decomposable measure for each
$c\in\overline{\mathbb{R}}_+$; \item[(ii)] the set function
$\zeta:=\tau(\gamma^1,\gamma^2)$ is a $\tau$-decomposable measure;
\item[(iii)] the set function
$\lambda:=\vartheta(\gamma^1,\gamma^2)$ is a $\tau$-decomposable
submeasure if and only if $\vartheta\gg\tau$.
\end{itemize}
\end{theorem}

\proof Consider $E, F\in\Sigma$ such that $E\cap F=\emptyset$. It
is easy to verify that $\gamma_\emptyset = \zeta_\emptyset =
\varepsilon_0$.

(i) Immediately, by distributivity of $\tau$ we get
$$\gamma_{E\cup F} = c\odot \gamma^1_{E\cup F} = c\odot \tau(\gamma^1_E,\gamma^1_F) =
\tau(c\odot \gamma^1_E, c\odot \gamma^1_F) =
\tau(\gamma_E,\gamma_F).$$

(ii) It follows from associativity of triangle functions that
\begin{align*}
\zeta_{E\cup F}& = \tau(\gamma^1_{E\cup F}, \gamma^2_{E\cup F})=
\tau(\tau(\gamma^1_E, \gamma^1_F), \gamma^2_{E\cup F})=
\tau(\gamma^1_E, \tau(\gamma^2_{E\cup F},
\gamma^1_F))=\tau(\gamma^1_E, \tau(\tau(\gamma^2_E, \gamma^2_F),
\gamma^1_F))\\
&= \tau(\gamma^1_E, \tau(\gamma^2_E,\tau(\gamma^2_F,
\gamma^1_F)))=\tau(\tau(\gamma^1_E,\gamma^2_E), \tau(\gamma^1_F,
\gamma^2_F)))= \tau(\zeta_E, \zeta_F).\end{align*}

(iii) If $\lambda:=\vartheta(\gamma^1,\gamma^2)$ is a
$\tau$-decomposable submeasure, then $\lambda_{E\cup
F}\geq\tau(\lambda_E,\lambda_F)$, i.e.,
$$\vartheta(\gamma_{E\cup F}^1,\gamma_{E\cup
F}^2)\geq\tau(\vartheta(\gamma_E^1,\gamma_E^2),\vartheta(\gamma_F^1,\gamma_F^2)).$$
Since
$$\vartheta(\gamma_{E\cup F}^1,\gamma_{E\cup
F}^2)=\vartheta(\tau(\gamma_E^1,\gamma_F^1),\tau(\gamma_E^2,\gamma_F^2)),
$$
therefore
\begin{equation}\label{dominance1}
\vartheta(\tau(\gamma_E^1,\gamma_F^1),\tau(\gamma_E^2,\gamma_F^2))\geq
\tau(\vartheta(\gamma_E^1,\gamma_E^2),\vartheta(\gamma_F^1,\gamma_F^2)),
\end{equation}
which means that $\vartheta\gg\tau$.

On the other hand, if $\vartheta\gg\tau$, then the
inequality~(\ref{dominance1}) holds and the equalities
\begin{align*}
\vartheta(\tau(\gamma_E^1,\gamma_F^1),\tau(\gamma_E^2,\gamma_F^2))
& = \vartheta(\gamma_{E\cup F}^1,\gamma_{E\cup
F}^2)=\lambda_{E\cup F}, \\
\tau(\vartheta(\gamma_E^1,\gamma_E^2),\vartheta(\gamma_F^1,\gamma_F^2))
& = \tau(\lambda_E,\lambda_F)
\end{align*} imply the result.\qed
\vskip 5pt

As we can see, a $\tau$-sum of two $\tau$-decomposable measures
produces a $\tau$-decomposable measure. The same is true for
(antimonotone) $\tau$-decomposable submeasures. Replacing a
triangle function by a more general aggregation operator (on
$\Delta^+$) may also produce a desired set function.

Recall that an ($n$-ary) \textit{aggregation operator} on
$\Delta^+$ is a mapping $\alpha:(\Delta^+)^n\rightarrow\Delta^+$
which is non-decreasing in each place with the boundary condition
$\alpha(\varepsilon_0,\dots,\varepsilon_0)=\varepsilon_0$. Also,
for domination between aggregation operators we refer
to~\cite{SMB}: we say that an ($n$-ary) aggregation operator
$\alpha$ \textit{dominates} an ($m$-ary) aggregation operator
$\beta$, we write $\alpha\gg\beta$, if for all
$H_{i,j}\in\Delta^+$ with $i\in\{1,\dots,m\}$ and
$j\in\{1,\dots,n\}$, the following property holds:
$$\alpha\bigl(\beta(H_{1,1}, \dots, H_{m,1}), \dots, \beta(H_{1,n}, \dots H_{m,n})\bigr) \geq
\beta\bigl(\alpha(H_{1,1},\dots, H_{1,n}), \dots,
\alpha(H_{m,1},\dots, H_{m,n})\bigr).$$

\begin{theorem}\label{dominates}
Let $\tau_i$ be triangle functions on $\Delta^+$ and
$\gamma^i:\Sigma\to\Delta^+$ be $\tau_i$-decomposable submeasures
for $i=1,\dots,n$, $n\in\mathbb{N}$. If $\tau$ is a triangle
function on $\Delta^+$ such that $\tau\leq \tau_i$ for each
$i=1,\dots, n$, and $\alpha$ is an ($n$-ary) aggregation operator
on $\Delta^+$ such that $\alpha\gg \tau$, then $\gamma:=
\alpha(\gamma^1, \cdots \gamma^n)$ is a $\tau$-decomposable
submeasure on $\Sigma$.
\end{theorem}

\proof Immediately, for each $E,F\in\Sigma$ we have
\begin{align*}
\gamma_{E\cup F}& = \alpha(\gamma^1_{E\cup F}, \dots,
\gamma^n_{E\cup F})\geq \alpha(\tau _1(\gamma^1_{E},
\gamma^1_{F}), \dots, \tau_n(\gamma^n_{E}, \gamma^n_{F}))\\ &\geq
\alpha(\tau(\gamma^1_{E}, \gamma^1_{F}), \dots, \tau(\gamma^n_{E},
\gamma^n_{F}))\geq \tau(\alpha(\gamma^1_{E},\cdots,\gamma^n_{E}),
\alpha(\gamma^1_{F},\cdots,\gamma^n_{F}))
\\ & =\tau(\gamma_E,\gamma_F).
\end{align*}Finally, by boundary condition we get $\gamma_\emptyset = \alpha(\gamma^1_\emptyset,\dots,\gamma^n_\emptyset)
=\alpha(\varepsilon_0,\dots,\varepsilon_0)=\varepsilon_0$. \qed
\vskip 5pt

\begin{remark}\rm
Especially, it is well-known that the arithmetic mean $AM$ on
$[0,1]$ dominates $W$, cf.~\cite{PT}, thus $\Pi_{AM}$ dominates
$\tau_{L,W}$ for an arbitrary $L\in\mathcal{L}$,
cf.~\cite{SamSemII}. So, $\Pi_{AM}$-aggregation of
$\tau_{L,W}$-decomposable submeasures, i.e., a set function
$\gamma:=\Pi_{AM}(\gamma^1,\dots,\gamma^n)$, is again a
$\tau_{L,W}$-decomposable submeasure.
\end{remark}

The following construction is based on the observation taken from
a construction of product PM-spaces for the case of finite
products, cf.~\cite{SamSemII}. For more details about
pseudo-metrics generated by probabilistic-valued decomposable
(sub)measures see below.

\begin{theorem}
Let $\tau_i$ be a (finite) family of triangle functions on
$\Delta^+$ and $\gamma^i:\Sigma_i\to\Delta^+$ be
$\tau_i$-decomposable antimonotone submeasures for $i=1,\dots\,n$,
$n\in\mathbb{N}$. Let $(\Sigma_i,\rho^i,\tau_i)$ be a family of
PpM-spaces with $\rho^i_{E_i,F_i}:=\gamma^i_{E_i\triangle F_i}$,
and $\alpha$ be an ($n$-ary) aggregation operator on $\Delta^+$.
If there exists a triangle function $\tau$ on $\Delta^+$ such that
$\alpha\gg\tau$ and that $\tau\leq\tau_i$ for every $i=1,\dots,n$,
then the triple $(\Sigma, \rho, \tau)$ is a PpM-space, where
$\rho$ is defined on $\Sigma:=\prod_{i=1}^n\Sigma_i$ by
$$\rho_{E,F}:=\alpha(\rho^1_{E_1,F_1},\dots,\rho^n_{E_n,F_n}).$$
\end{theorem}

\proof We prove that $\rho$ is a pseudo-metric on $\Sigma$ w.r.t.
$\tau$. Clearly, $\rho_{E,E}=\varepsilon_0$ for each $E\in\Sigma$
and $\rho_{E,F}=\rho_{F,E}$ for each $E,F\in\Sigma$. Taking
$E,F,G\in\Sigma$ we obtain
\begin{align*}
\tau(\rho_{E,F},\rho_{F,G}) & =
\tau(\alpha(\rho^1_{E_1,F_1},\dots,\rho^n_{E_n,F_n}),\alpha(\rho^1_{F_1,G_1},\dots,\rho^n_{F_n,G_n}))\\
& \leq
\alpha(\tau(\rho^1_{E_1,F_1},\rho^1_{F_1,G_1}),\dots,\tau(\rho^n_{E_n,F_n},\rho^n_{F_n,G_n}))\\
& \leq
\alpha(\tau_1(\rho^1_{E_1,F_1},\rho^1_{F_1,G_1}),\dots,\tau_n(\rho^n_{E_n,F_n},\rho^n_{F_n,G_n}))\\
& \leq
\alpha(\rho^1_{E_1,G_1},\dots,\rho^n_{E_n,G_n})=\rho_{E,G},\
\end{align*}which completes the proof.\qed \vskip 5pt


We now extend results from~\cite{HalHutMes} to a general case of
$\tau$-decomposable (sub)measures instead of
$\tau_{L,T}$-submeasures considered therein.

\begin{theorem}\label{pseudometricNA}
Let $\tau$ be a triangle function on $\Delta^+$ and $\gamma:
\Sigma\to\Delta^+$ be a $\tau$-decomposable antimonotone
submeasure. Define the mapping $\rho:
\Sigma\times\Sigma\to\Delta^+$ by the formula
$\rho_{E,F}:=\gamma_{E\triangle F}$ with $E, F\in\Sigma$. Then the
triple $\left(\Sigma, \rho, \tau\right)$ is a PpM-space (under
$\tau$).
\end{theorem}

\proof Clearly, $\rho_{E,E}=\gamma_{E\triangle E}=\varepsilon_0$
for each $E\in\Sigma$. Similarly, symmetry is trivial. It is
enough to prove the triangle inequality. By
Theorem~\ref{thmantimonotone} consider arbitrary sets $E, F,
G\in\Sigma$. Since
$$E\triangle G\subset (E\triangle F)\cup [(G\setminus(E\cup
F))\cup ((E\cap F)\setminus G)],
$$ and
$$(G\backslash(E\cup F))\cup((E\cap F)\backslash G)\subset
F\triangle G,$$ then antimonotonicity of $\gamma$ yields
\begin{align}\label{triangleineq}
\rho_{E,G} =\gamma_{E\triangle G}& \geq \gamma_{(E\triangle F)\cup
[(G\setminus(E\cup F))\cup ((E\cap F)\setminus G)]}\geq
\tau(\gamma_{E\triangle
F},\gamma_{(G\setminus(E\cup F))\cup ((E\cap F)\setminus G)})\nonumber\\
& \geq\tau(\gamma_{E\triangle F},\gamma_{F\triangle
G})=\tau(\rho_{E,F},\rho_{F,G}),\ \end{align} which means that
$\rho$ is a pseudo-metric with respect to the triangle function
$\tau$.\qed \vskip 5pt


\begin{remark}\rm
Since $E\triangle G = (E\triangle F)\triangle (F\triangle G)$,
then it is easy to observe that $\rho$ is translation invariant,
i.e., $\rho_{E,G} = \gamma_{E\triangle G} = \gamma_{(E\triangle
F)\triangle (F\triangle G)} = \rho_{E\triangle F, F\triangle G}$.
\end{remark}

For a fixed triangle function $\tau$ denote by
$\Gamma_{\tau}(\Sigma, \gamma)$ \textit{the set of all
pseudo-metrics $\rho$ generated by $\tau$-decomposable
antimonotone submeasures} $\gamma$ on $\Sigma$, i.e.,
$\rho_{E,F}=\gamma_{E\triangle F}$ with $E,F\in\Sigma$. Define a
relation $\preceq$ on $\Gamma_{\tau}(\Sigma, \gamma)$ as follows
$$\rho \preceq \varrho \Leftrightarrow \rho_{E,F} \geq \varrho_{E,F}\,\,\,\textrm{for each}\,\,\, E,F\in\Sigma.$$ Clearly, $\preceq$ is a
partial order on $\Gamma_{\tau}(\Sigma, \gamma)$ and
$\nu_{E,F}=\varepsilon_0$ for each $E,F\in\Sigma$ is an element
such that $\nu\preceq \rho$ for every
$\rho\in\Gamma_{\tau}(\Sigma, \gamma)$.

Let $\vartheta$ be a triangle function on $\Delta^+$ such that
$\vartheta\gg \tau$. Define a binary operation
$\oplus_{\vartheta}$ on $\Gamma_{\tau}(\Sigma, \gamma)$ such that
$$(\rho\oplus_{\vartheta}\varrho)_{E,F} := \vartheta(\rho_{E,F},
\varrho_{E,F}), \quad E,F\in \Sigma,$$ for all $\rho, \varrho \in
\Gamma_{\tau}(\Sigma, \gamma)$. We show that
$\rho\oplus_{\vartheta}\varrho \in \Gamma_{\tau}(\Sigma, \gamma)$.
Indeed, for an arbitrary $E\in\Sigma$ we have
$$(\rho\oplus_{\vartheta}\varrho)_{E,E} = \vartheta(\rho_{E,E},
\varrho_{E,E}) = \vartheta(\varepsilon_0, \varepsilon_0) =
\varepsilon_0$$ and the symmetry
$(\rho\oplus_{\vartheta}\varrho)_{E,F}=(\rho\oplus_{\vartheta}\varrho)_{F,E}$
is obvious for each $E,F\in\Sigma$. Also, by triangle
inequality~(\ref{triangleineq}) and domination of $\vartheta$ over
$\tau$ we have {\setlength\arraycolsep{2pt}
\begin{eqnarray*}
(\rho\oplus_{\vartheta}\varrho)_{E,G} & = & \vartheta(\rho_{E,G},
\varrho_{E,G}) \geq \vartheta(\tau(\rho_{E,F}, \rho_{F,G}),
\tau(\varrho_{E,F}, \varrho_{F,G})) \\ & \geq &
\tau(\vartheta(\rho_{E,F}, \varrho_{E,F}), \vartheta(\rho_{F,G},
\varrho_{F,G})) \\ & = &
\tau((\rho\oplus_{\vartheta}\varrho)_{E,F},
(\rho\oplus_{\vartheta}\varrho)_{F,G}),
\end{eqnarray*}}where $E,F,G\in\Sigma$. Thus $\rho\oplus_{\vartheta}\varrho$ is a pseudo-metric on $\Sigma$.
Note that for all $\rho\in\Gamma_{\tau}(\Sigma, \gamma)$ and
$\vartheta\gg \tau$ it holds
$$(\nu\oplus_{\vartheta}\rho)_{E,F} = \vartheta(\nu_{E,F},
\rho_{E,F}) = \vartheta(\varepsilon_0,\rho_{E,F}) = \rho_{E,F}$$
for each $E,F\in\Sigma$. The operation $\oplus_{\vartheta}$ is
clearly commutative and associative. Thus we have the following
result.


\begin{proposition}\label{propositiontriple}
Let $\vartheta\gg \tau$. The triple $(\Gamma_{\tau}(\Sigma,
\gamma), \oplus_{\vartheta}, \preceq)$ is a partially ordered
commutative semigroup with the neutral element $\nu$.
\end{proposition}

A relationship between the relation $\preceq$ and the operation
$\oplus_{\vartheta}$ is described in the following theorem.

\begin{theorem}\label{thmoplusprec}
Let $\vartheta\gg \tau$. Then for all $\rho, \varrho, \sigma \in
\Gamma_{\tau}(\Sigma, \gamma)$ the relationship
$\rho\oplus_{\vartheta} \sigma \preceq \varrho\oplus_{\vartheta}
\sigma$ holds whenever $\rho \preceq \varrho$.
\end{theorem}

\proof The relation $\rho \preceq \varrho$ means that
$\rho_{E,F}\geq \varrho_{E,F}$ for each $E,F\in\Sigma$. From
monotonicity of $\vartheta$ we have $\vartheta(\rho_{E,F},
\sigma_{E,F})\geq \vartheta(\varrho_{E,F}, \sigma_{E,F})$ and the
proof is complete. \qed \vskip 5pt

Recall that $\Pi_M$ (the pointwise defined operation on $\Delta^+$
given by~(\ref{tautop}) with the minimum t-norm $M$) dominates
each triangle function $\tau$, cf.~\cite[Proposition
6.2]{SamSemII}. 
Also, a \textit{semilattice} is an idempotent commutative
semigroup. A semilattice is \textit{bounded} if it includes the
neutral element.

\begin{theorem}
The ordered pair $(\Gamma_{\tau}(\Sigma, \gamma), \oplus_{\Pi_M})$
is a bounded semilattice with the properties
\begin{itemize}
\item[(i)] $\rho \preceq \varrho$ if and only if $\rho
\oplus_{\Pi_M} \varrho = \varrho$; \item[(ii)]
$(\sigma\oplus_{\vartheta} \rho) \oplus_{\Pi_M}
(\sigma\oplus_{\vartheta} \varrho) \preceq \sigma
\oplus_{\vartheta}(\rho\oplus_{\Pi_M} \varrho)$
\end{itemize}\noindent for each $\rho, \varrho, \sigma \in \Gamma_{\tau}(\Sigma, \gamma)$ provided that $\vartheta\gg \tau$.
\end{theorem}

\proof Since $\nu\in \Gamma_{\tau}(\Sigma, \gamma)$ is the neutral
element, then $(\Gamma_{\tau}(\Sigma, \gamma), \oplus_{\Pi_M})$ is
a bounded semilattice.

(i) If $\rho\preceq\varrho$, then $\rho_{E,F}\geq \varrho_{E,F}$
for each $E,F\in\Sigma$, thus $\rho \oplus_{\Pi_M} \varrho =
\varrho$. The opposite statement is obvious.

(ii) For the proof of second part it is enough to observe that
$\rho \preceq \rho \oplus_{\Pi_M} \varrho$ and $\varrho \preceq
\rho\oplus_{\Pi_M}\varrho$. By Theorem~\ref{thmoplusprec} we have
$\sigma\oplus_{\vartheta} \rho \preceq
\sigma\oplus_{\vartheta}(\rho\oplus_{\Pi_M}\varrho)$ and
$\sigma\oplus_{\vartheta} \varrho \preceq
\sigma\oplus_{\vartheta}(\rho\oplus_{\Pi_M}\varrho)$. Since
$(\Gamma_{\tau}(\Sigma, \gamma), \oplus_{\Pi_M})$ is a
semilattice, the property~(ii) follows. \qed \vskip 5pt

\section{Decomposable measures and the probabilistic Hausdorff
distance}\label{sec_hausdorff}

As it was already mentioned in this paper, Shen introduced
in~\cite{Shen} a class of probabilistic decomposable measures
w.r.t. the triangle function $\Pi_\top$ given by~(\ref{tautop}).
Now we extend his certain results related to the probabilistic
Hausdorff distance. For that reason we recall necessary notions
from~\cite{Shen}.\vskip 5pt

Let $(\Omega,\mathcal{F},\tau)$ be a PM-space.

\noindent (i) The \textit{probabilistic diameter} of a non-empty
subset $E$ of $\Omega$ is a mapping $D_E: [0,+\infty)\to[0,1]$
defined by
$$D_E(t):=\sup_{s<t} \inf_{p,q\in E} H_{p,q}(s).$$

\noindent (ii) A set $E$ is said to be \textit{probabilistic
bounded} if $\sup_t D_E(t)=1$. The collection of all probabilistic
bounded subsets of $\Omega$ will be denoted by $PB(\Omega)$.

\noindent (iii) Let $p\in\Omega$ and $F\in PB(\Omega)$. The
\textit{probabilistic distance} from $p$ to $F$ is defined as
$$d_{p,F}(t) :=
\begin{cases}
0, & t=0, \\
\sup\limits_{s<t} \sup\limits_{q\in F} H_{p,q}(s), & t\in
\mathbb{R}_+,
\end{cases}$$ with the convention
$d_{p,\emptyset}=1-\varepsilon_0$.

\noindent (iv) Given $E,F\in PB(\Omega)$, the
\textit{probabilistic distance} from $E$ to $F$ is defined as
$$d_{E,F}(t) :=
\begin{cases}
0, & t=0, \\
\sup\limits_{s<t} \inf\limits_{p\in E} \sup\limits_{q\in F}
H_{p,q}(s), & t\in \mathbb{R}_+,
\end{cases}$$ with the convention
$d_{\emptyset,F}=\varepsilon_0$.

\begin{definition}\rm
Let $(\Omega,\mathcal{F},\tau)$ be a PM-space and $E,F\in
PB(\Omega)$. The probabilistic Hausdorff distance between $E$ and
$F$ is a mapping $\mathcal{H}_{E,F}: [0,+\infty)\to [0,1]$ defined
by $$\mathcal{H}_{E,F}(t):=\begin{cases}
0, & t=0, \\
\sup\limits_{s<t}\, M\left(\inf\limits_{p\in E} \sup\limits_{q\in
F} H_{p,q}(s), \inf\limits_{q\in F} \sup\limits_{p\in E}
H_{p,q}(s)\right), & t\in \mathbb{R}_+.
\end{cases}$$
\end{definition}

\noindent From now on let $\Sigma:= P(\Omega)$ be the power set of
$\Omega$. Using the probabilistic Hausdorff distance we may define
a mapping $\Lambda^\mathcal{H}: \Sigma\to\Delta^+$ by
$$\Lambda^\mathcal{H}_X := \mathcal{H}_{E^c, \Omega},$$ where $E^c$
is the complement of a set $E\in \Sigma$.

\begin{theorem}
Let $(\Omega,\mathcal{F},\tau)$ be a PM-space. Then
$\Lambda^\mathcal{H}$ is an antimonotone set function with
$\Lambda^\mathcal{H}_\emptyset=\varepsilon_0$ and
$$\Lambda^\mathcal{H}_{E\cup F} \leq \Pi_M\left(\Lambda^\mathcal{H}_E, \Lambda^\mathcal{H}_F\right), \quad E,F\in \Sigma.$$
\end{theorem}

\proof 
From definitions above it is obvious that
$\Lambda^\mathcal{H}_\emptyset = \varepsilon_0$. Also, for $E,F\in
\Sigma$ such that $E\subseteq F$ we have $F^c\subseteq E^c$, and
therefore
$$\Lambda^\mathcal{H}_{E}(t) = \mathcal{H}_{E^c,\Omega}(t) = \sup_{s<t} \inf_{q\in\Omega} \sup_{p\in E^c} H_{p,q}(s) \geq
\sup_{s<t} \inf_{q\in\Omega} \sup_{p\in F^c} H_{p,q}(s) =
\mathcal{H}_{F^c,\Omega}(t) = \Lambda^\mathcal{H}_{F}(t)$$ for
each $t\geq 0$, which proves antimonotonicity of
$\Lambda^\mathcal{H}$.

Finally, the inclusions $E\subseteq E\cup F$, $F\subseteq E\cup F$
and antimonotonicity of $\Lambda^\mathcal{H}$ yields
$$\Lambda^\mathcal{H}_{E\cup F}(t)\leq M\left(\Lambda^\mathcal{H}_{E}(t), \Lambda^\mathcal{H}_{F}(t)\right) =
\Pi_M\left(\Lambda^\mathcal{H}_{E}(t),
\Lambda^\mathcal{H}_{F}(t)\right)$$ for each $t\geq 0$, i.e.,
$\Lambda^\mathcal{H}_{E\cup F} \leq
\Pi_M\left(\Lambda^\mathcal{H}_E, \Lambda^\mathcal{H}_F\right)$.
\qed

\begin{remark}\rm
Equivalently, we may say that $\Lambda^\mathcal{H}$ is a
\textit{probabilistic-valued antimonotone $\Pi_M$-decomposable
supermeasure on $\Sigma$} ($\Pi_M$-probabilistic submeasure in the
terminology of Shen~\cite{Shen}). Under certain condition on the
domain of $\Lambda^\mathcal{H}$ we will be able to show that it is
a $\tau$-decomposable measure w.r.t. an arbitrary triangle
function $\tau$ (thus, also for $\Pi_M$), see
Theorem~\ref{thmrestriction} below.
\end{remark}

\begin{definition}\rm
A set $E\in\Sigma$ is said to be
\textit{$\Lambda^\mathcal{H}$-measurable w.r.t. a triangle
function} $\tau$, if for each $G\in\Sigma$ it holds
$$\Lambda^\mathcal{H}_G = \tau\left(\Lambda^\mathcal{H}_{G\cap E},
\Lambda^\mathcal{H}_{G\cap E^c}\right).$$
\end{definition}

\begin{theorem}
Let $(\Omega,\mathcal{F},\tau)$ be a PM-space. Then each class
$$\mathbb{S}_\tau=\{E\in\Sigma;\,\,
E\,\,\textrm{is}\,\,\Lambda^\mathcal{H}\textrm{-measurable}\,\,\textrm{w.r.t.}\,\,\tau\}$$
is an algebra.
\end{theorem}

\proof Since $\mathcal{H}_{\emptyset,\Omega} =
\mathcal{H}_{\Omega,\Omega}=\varepsilon_0$, then $\emptyset,
\Omega\in\mathbb{S}_\tau$. Moreover, from definition it follows
that if $E\in\mathbb{S}_\tau$, then $E^c\in\mathbb{S}_\tau$ (and
vice versa). Thus, we only have to prove that $\mathbb{S}_\tau$ is
closed under the formation of union.

Let $E,F\in\mathbb{S}_\tau$ and consider an arbitrary
$G\in\Sigma$. Since $F\in\mathbb{S}_\tau$, then
$$\Lambda^\mathcal{H}_{G\cap (E\cup F)} =
\tau\left(\Lambda^\mathcal{H}_{G\cap(E\cup F)\cap F},
\Lambda^\mathcal{H}_{G\cap(E\cup F)\cap F^c}\right) =
\tau\left(\Lambda^\mathcal{H}_{G\cap F},
\Lambda^\mathcal{H}_{G\cap E\cap F^c}\right),$$ and therefore we
get
\begin{align*}
\tau\left(\Lambda^\mathcal{H}_{G\cap(E\cup F)},
\Lambda^\mathcal{H}_{G\cap(E\cup F)^c}\right) & =
\tau\left(\tau\left(\Lambda^\mathcal{H}_{G\cap F},
\Lambda^\mathcal{H}_{G\cap F^c\cap E}\right),
\Lambda^\mathcal{H}_{G\cap(E\cup F)^c}\right) \\ & =
\tau\left(\tau\left(\Lambda^\mathcal{H}_{G\cap F},
\Lambda^\mathcal{H}_{G\cap F^c\cap E}\right),
\Lambda^\mathcal{H}_{G\cap F^c\cap E^c}\right) \\ & =
\tau\left(\Lambda^\mathcal{H}_{G\cap F},
\tau\left(\Lambda^\mathcal{H}_{G\cap F^c\cap E},
\Lambda^\mathcal{H}_{G\cap F^c\cap E^c}\right)\right) \\ & =
\tau\left(\Lambda^\mathcal{H}_{G\cap F},
\Lambda^\mathcal{H}_{G\cap F^c}\right) \\ & =
\Lambda^\mathcal{H}_G.
\end{align*}Thus, the set $E\cup F$ is $\Lambda^\mathcal{H}$-measurable w.r.t.
$\tau$, i.e., $E\cup F\in\mathbb{S}_\tau$. \qed

\begin{theorem}\label{thmrestriction}
Let $(\Omega,\mathcal{F},\tau)$ be a PM-space. The restriction of
the set function $\Lambda^\mathcal{H}$ to $\mathbb{S}_\tau$ is a
$\tau$-decomposable measure.
\end{theorem}

\proof Let $E,F\in\mathbb{S}_\tau$ be disjoint. Since $E\cup F$ is
$\Lambda^\mathcal{H}$-measurable w.r.t. $\tau$, then
$$\Lambda^\mathcal{H}_{E\cup F} = \tau\left(\Lambda^\mathcal{H}_{(E\cup F)\cap E}, \Lambda^\mathcal{H}_{(E\cup F)\cap E^c}\right)
= \tau\left(\Lambda^\mathcal{H}_{E},
\Lambda^\mathcal{H}_{F}\right),$$ i.e., $\Lambda^\mathcal{H}$ is a
$\tau$-decomposable measure on $\mathbb{S}_\tau$. \qed

\section*{Concluding remarks}

In this paper we have introduced a generalized decomposable
(sub)measure notion related to triangle functions in probabilistic
metric spaces which covers our previously introduced notions of
(sub)measures as well as certain recent notions from the
literature. We have given a characterization of decomposable
measures and discussed construction methods providing new
decomposable (sub)measures from the given ones. In connection with
this "aggregation" the following problem arises:

\begin{problem}
Characterize all the mappings which preserve the class of
$\tau$-decomposable measures.
\end{problem}

Further, we have extended Shen's results related to
probabilistic-valued decomposable measures and the probabilistic
Hausdorff distance for an arbitrary triangle function. We have
proved that the class of all measurable sets forms an algebra on
its power set. Thus, the second problem repeats the Shen's
question in~\cite{Shen}:

\begin{problem}
Is the class $\mathbb{S}_\tau$ of all
$\Lambda^\mathcal{H}$-measurable sets w.r.t. a triangle function
$\tau$ a $\sigma$-algebra? Does this property depend on a choice
of triangle function $\tau$?
\end{problem}

As it is well-known, triangle functions have been used and
discussed almost exclusively in the theory of probabilistic metric
and probabilistic normed spaces. Particular examples of triangle
functions also appear in, e.g., the treatment of fuzzy numbers or
in information theory (compare also~\cite{PapSta}). As it was
described in the introduction the probabilistic-valued
(sub)measures w.r.t. triangle function $\tau_M$ may be naturally
interpreted in the context of fuzzy sets and fuzzy numbers. Thus,
one can be interested in the following:

\begin{problem}
Is there room in these theories (both PM-spaces as well as fuzzy
sets) for $\tau$-decomposable (sub)measures w.r.t. triangle
functions of other type? Which interpretation do this objects
have?
\end{problem}

Furthermore, the idea of probabilistic integral introduced
in~\cite{Lpspaces} motivates to develop an integral w.r.t.
probabilistic-valued decomposable (sub)measures for an arbitrary
triangle function. This is already done in our recent
paper~\cite{HalHut}. This approach may be understood as a
modification of Aumann integral, or, of Choquet-type integral
based on interval-valued measures as discussed e.g.
in~\cite{Jang}.

\section*{Acknowledgement}

We are grateful to the referees for their suggestions which
essentially improved previous version of the manuscript. We
acknowledge a partial support of grants VEGA 1/0171/12,
VVGS-PF-2013-115 and VVGS-2013-121.


\vspace{5mm}

\noindent \small{\textsc{Lenka Hal\v cinov\'a, Ondrej Hutn\'ik,
Jana Moln\'arov\'a} \newline Institute of Mathematics, Faculty of
Science, Pavol Jozef \v Saf\'arik University in Ko\v sice,
\newline {\it Current address:} Jesenn\'a 5, SK 040~01 Ko\v sice,
Slovakia,
\newline {\it E-mail addresses:} lenka.halcinova@student.upjs.sk \newline
\phantom{{\it E-mail addresses:}} ondrej.hutnik@upjs.sk \newline
\phantom{{\it E-mail addresses:}} jana.molnarova@student.upjs.sk}


\begin{thebibliography}{1}

\bibitem{dBM}
\textsc{de Baets, B., Mesiar, R.}: $\mathcal{T}$-partitions.
\textit{Fuzzy Sets and Systems}~\textbf{97} (1998), 211--223.

\bibitem{Lpspaces}
\textsc{Bahrami, F., Mohammadbaghban, M.}: Probabilistic $L^p$
spaces. \textit{J. Math. Anal. Appl.}~\textbf{402}(2) (2013),
505--517.

\bibitem{DKMP}
\textsc{Dubois, D. Kerre, E.~E., Mesiar, R. Prade, H.}: Fuzzy
interval analysis. In: \textit{Fundamentals of Fuzzy Sets}, D.
Dubois and H. Prade (Eds.), Kluwer, Dordrecht, 2000, pp. 483--584.

\bibitem{GS}
\textsc{Gilio, A, Sanfilippo, G.}: Quasi conjunction, quasi
disjunction, t-norms and t-conorms: probabilistic aspects.
\textit{Inf. Sci.}~\textbf{245} (2013), 146--167.

\bibitem{GPMM}
\textsc{Grabisch, M., Pap, E., Mesiar, R., Marichal, J.-L.}:
\textit{Aggregation Functions.} Cambridge Univ. Press, New York,
2009.

\bibitem{GCK}
\textsc{Grzymala-Busse, J. W., Clark, P. G., Kuehnhausen, M.}:
Generalized probabilistic approximations of incomplete data.
\textit{Internat. J. Approx. Reason.}~\textbf{55} (2014),
180--196.

\bibitem{HalHut}
\textsc{Hal\v cinov\'a, L., Hutn\'ik}: An integral with respect to
probabilistic-valued decomposable measures. (submitted).

\bibitem{HalHutMes2}
\textsc{Hal\v cinov\'a, L., Hutn\'ik, O., Mesiar, R.}: On distance
distribution functions-valued submeasures related to aggregation
functions. \textit{Fuzzy Sets and Systems}~{\bf 194}(1) (2012),
15--30.

\bibitem{HalHutMes}
\textsc{Hal\v cinov\'a, L., Hutn\'ik, O., Mesiar, R.}: On some
classes of distance distribution functions-valued submeasures.
\textit{Nonlinear Anal.}~{\bf 74}(5) (2011), 1545--1554.

\bibitem{Hunter}
\textsc{Hunter, A.}: A probabilistic approach to modelling
uncertain logical arguments. \textit{Internat. J. Approx.
Reason.}~\textbf{54} (2013), 47--81.

\bibitem{HutMes}
\textsc{Hutn\'ik, O., Mesiar, R.}: On a certain class of
submeasures based on triangular norms. \textit{Internat. J.
Uncertain. Fuzziness Knowledge-Based Systems}~\textbf{17} (2009),
297--316.

\bibitem{Jang}
\textsc{Jang, L.-Ch.}: A note on convergence properties of
interval-valued capacity functionals and Choquet integrals.
\textit{Inf. Sci.}~\textbf{183}(1) (2012), 151--158.

\bibitem{KdF74}
\textsc{Kamp\'e de F\'eriet, J.-M.}: La th\'eorie general de
l'information et la mesure subjective de l'information. In:
\textit{Lecture Notes in Math.}~{\bf 398}, Springer-Verlag,
Heidelberg 1974, pp. 1--35.

\bibitem{KMP}
\textsc{Klement, E.~P., Mesiar, R., Pap, E.}: \textit{Triangular
Norms.} Trends in Logic, Studia Logica Library {\bf 8}, Kluwer
Academic Publishers, 2000.

\bibitem{Menger}
\textsc{Menger, K.}: Statistical metrics. \textit{Proc. Nat. Acad.
Sci. U.S.A.}~{\bf 28} (1942), 535--537.

\bibitem{MB}
\textsc{Moore, R. E., Bierbaum, F.}: \textit{Methods and
Applications of Interval Analysis}. Studies in Applied and
Numerical Mathematics, SIAM, Philadelphia, 1979.

\bibitem{PapSta}
\textsc{Pap, E., \v{S}tajner, I.}: Generalized pseudo-convolution
in the theory of probabilistic metric spaces, information, fuzzy
numbers, optimization, system theory. \textit{Fuzzy Sets and
Systems}~{\bf 102} (1999), 393--415.

\bibitem{PT}
\textsc{Pradera, A., Trillas, E.}: A note on pseudometrics
aggregation, {\it Int. J. Gen. Syst.}~{\bf 31} (2002), 41--51.

\bibitem{SMB}
\textsc{Saminger, S., Mesiar, R., Bodenhofer, U.}: Domination of
aggregation operators and preservation of transitivity.
\textit{Internat. J. Uncertain. Fuzziness Knowledge-Based
Systems}~\textbf{10} (2002), 11--35.


\bibitem{SamSem}
\textsc{Saminger-Platz, S. Sempi, C.}: A primer on triangle
functions I., {\it Aequationes Math.}~{\bf 76}(3) (2008),
201--240.

\bibitem{SamSemII}
\textsc{Saminger-Platz, S. Sempi, C.}: A primer on triangle
functions II., {\it Aequationes Math.}~{\bf 80}(3) (2010),
239--268.

\bibitem{SS}
\textsc{Schweizer, B., Sklar, A.}: \textit{Probabilistic Metric
Spaces.} North-Holland Publishing, New York, 1983.

\bibitem{Shen}
\textsc{Shen, Y.}: On probabilistic Hausdorff distance and a class
of probabilistic decomposable measures. \textit{Inf.
Sci.}~\textbf{263} (2014), 126--140.

\bibitem{SM}
\textsc{Sugeno, M., Murofushi, T.}:  Pseudo-additive measures and
integrals. \textit{J. Math. Anal. Appl.}~\textbf{122} (1987),
197--222.

\bibitem{cinania}
\textsc{Yu, Z., Li, L. Wong, H.-S., You, J., Han, G., Gao, Y., Yu,
G.}: Probabilistic cluster structure ensemble. \textit{Inf. Sci.},
http://dx.doi.org/10.1016/j.ins.2014.01.030.
\end{thebibliography}
\end{document}